\documentclass[a4paper,12pt]{article}
\usepackage[latin1]{inputenc}
\usepackage{amsmath,amssymb,amsthm,cite,verbatim}
\usepackage[all]{xy}
\usepackage[bookmarks=true]{hyperref}
\usepackage{color}
\usepackage{graphicx}

\title{The automorphism group of the  modular curve $X_0^*(N)$ with square-free level}

\author{ \ Francesc Bars\footnote{First author is supported by MTM2016-75980-P} and Josep Gonz\'alez  \footnote{The  second author is partially supported by DGI grant  MTM2015-66180-R.}}

\newtheorem{prop}{Proposition}
\newtheorem{lema}{Lemma}
\newtheorem{teo}{Theorem}
\newtheorem{cor}{Corollary}
\newtheorem{rem}{Remark}

\theoremstyle{definition}

\theoremstyle{remark}

\numberwithin{equation}{section}

\newcommand{\Q}{\mathbb{Q}}

\newcommand{\Z}{\mathbb{Z}}

\newcommand{\F}{\mathbb{F}}

\newcommand{\PP}{\mathbf{P}}
\newcommand{\C}{\mathbb{C}}

\newcommand{\R}{\mathbb{R}}

\newcommand{\Gal}{\mathrm{Gal}}

\newcommand{\SL}{\operatorname{SL}}
\newcommand{\End}{\operatorname{End}}
\newcommand{\Aut}{\operatorname{Aut}}

\newcommand{\Frob}{\operatorname{Frob}}
\newcommand{\Jac}{\operatorname{Jac}}

\newcommand{\New}{\operatorname{New}}

\newcommand{\cO}{{\mathcal O}}

\newcommand{\cI}{{\mathcal I}}

\newcommand{\cS}{{\mathcal S}}
\newcommand{\cP}{{\mathcal P}}

\newcommand{\cL}{{\mathcal L}}

\date{}
\begin{document}
\maketitle

\begin{abstract}
\noindent
We determine the automorphism group of the modular curve $X_0^*(N)$, obtained as the quotient of the modular curve $X_0(N)$ by the group of its Atkin-Lehner involutions, for all square-free values of $N$.
\end{abstract}

\section{Introduction}
In \cite{KeMo}, Kenku and Momose  determined the automorphism group of all modular curves $X_0(N)$ with genus $>1$, except to $N=63$ that was solved by Elkies in \cite{elkies2}.  Subsequently, Harrison detected  a mistake  in these results concerning to  the curve $X_0(108)$. In  \cite{Har}, it is proved that this curve  has an extra involution that does  not arise from  the normalizer of $\Gamma_0(108)$ in $\SL_2(\R)$, as it happens for the curves $X_0(37)$ and $X_0(63)$.

For the modular curve $X_0^+(N)=X_0(N)/\langle w_N\rangle$, where  $w_N$ denotes the Fricke involution, Baker and Hasegawa determined the automorphism group when $N$ is a prime in \cite{BH}] and, later, in \cite{Go16} it was determined the automorphism group when $N$ is the square of a prime.

In this paper we focus our attention to the modular curves $X_0^*(N)=X_0(N)/B(N)$, where  $B(N)$ is the group of the Atkin-Lehner involutions
of the modular curve $X_0(N)$, when $N$ is square-free. The interest in these modular curves is due to their moduli interpretation. For a number field $K$, the non cuspidal $K$-rational points of $X_0^*(N)$ parametrize a class of
 $K$-curves. More precisely, these points parametrize elliptic curves  $E/\overline K$ having the property  that for every Galois conjugation $\sigma \in\Gal (\overline K/K)$ there is an isogeny between $E$ and $E^\sigma$ of degree an divisor of $N$. When  $X_0^*(N)$ has genus $<3$, these parametrizations are described in   \cite{GL} and \cite{BGX}.

We point out that for a square-free $N\neq 37$, $B(N)$ is the
automorphism group of $X_0(N)$ when its genus is  greater than one
and, thus,  a non trivial automorphism of $X_0^*(N)$ does not come
from the action of an automorphism of $X_0(N)$.
 {{More
generally, let us denote by $A(N)$ the subgroup of $\SL_2(\R)$  generated by $\Gamma_0(N)$ and $B(N)$. In \cite{LaMo}, it is proved that $A(N)$ is its normalizer in $\SL_2(\R)$. Hence,
 a non trivial automorphism of
$X_0^*(N)$ is exceptional  in the sense that does not come from a linear
fractional transformation on the complex upper-half plane.}}

In  \cite{BH}[Corollary 2.6], it is proved
that for a square-free integer $N$ such that  the genus of  $X_0^*(N)$ is $>1$,
 the group $\Aut X_0^*(N)$ is  elementary $2$-abelian and every automorphism of $X_0^*(N)$ is defined over $\Q$.
If  the curve $X_0^*(N)$ has a non trivial involution and the genus of the quotient curve is zero, then it is hyperelliptic.  When the genus of the quotient curve is one, then the curve $X_0^*(N)$   is bielliptic.

 In \cite{HaHa}, it is proved that $X_0^*(N)$ is hiperelliptic if, and only if, the curve  has genus two and, moreover,  all these   values of $N$   are determined.  In \cite{BaGon}, all  values of $N$ for which $X_0^*(N)$  is bielliptic  are determined  and, moreover,   for everyone of these curves  it is determined its automorphism  group.

\vskip 0.2cm
All these results are summarized in the next theorem.

\begin{teo}\label{main} Let $N>1$ be a square-free integer. Assume that the genus of the modular curve $X_0^*(N)$
is at least $2$. Then,
\begin{itemize}
\item  [\rm{(i)}] The modular curve $X_0^*(N)$ is hyperelliptic if, and only if, has genus two. In this case,  the  automorphism group has order $2$
if, and only if, $N$ is in the following table
$$
\begin{array}{r}
{ 67}, { 73}, {85},  { 93}, { 103},
 { 107},  {115},
   { 133}, { 134},  {
146},   154,  {161}, { 165},
{167},  { 170}, { 177},\\{
186}, {191}, {205}, {206},
{ 209}, { 213},  { 221}, { 230},
{ 266} , {285}, {286}, { 287}, {299},
 { 357}.\\
\end{array}
$$

\item[\rm{ (ii)}] The automorphism group of a modular curve $X_0^*(N)$ of genus $3$ is non trivial if, and only if, the curve is bielliptic.
\item [\rm{(iii)}] The modular curve $X_0^*(N)$ is bielliptic
if, and only if, $N$ is in the following table
$$ \begin{array}{c|r|}
\text{genus} & N\phantom{cccccccccccccccccccccccc}\\\hline
2& 106,122,129, 158,166,215,390\\ \hline
3&  178,183, 246, 249, 258, 290, 303,318, 430,455,510\\\hline
4 & 370\\ \hline
\end{array}
$$
 For these values of $N$, the automorphism group of
$X_0^*(N)$  has order $2$ when its genus is greater than two,
otherwise it is the Klein group.

\end{itemize}
\end{teo}
The goal of this article is to complete the values of $N$ such that the group  $\Aut (X_0^*(N))$ is non trivial and describe   this  for all these values.
Among the curves $X_0^*(N)$ with genus $>1$, there are exactly $37$ that are  hyperelliptic and $12$ that are  bielliptic and non hyperelliptic.
In \cite{BaGon}, it is proved that the  curve  $X_0^*(366$) of genus $4$ has  automorphism group of order $2$ and the genus of the quotient curve by the non trivial involution is $2$. Hence, it is reasonable to expect that there are  a few curves $X_0^*(N)$ of genus $>3$ that are not bielliptic and  its automorphism group is non trivial.
The main result of this paper, that is  presented in the following theorem, gives a precise  answer to this question.
\begin{teo}\label{main2} Let $N$ be a square-free integer  such that the curve $X_0^*(N)$ has genus greater than $3$ and is not bielliptic, i.e. $N\neq 370$. Then, the group  $\Aut (X_0^*(N))$ is not trivial, if and only if, $N= 366, 645$. In both cases, the order of
 this group is $2$ and the genus of the quotient curve by the non trivial involution is $2$.
\end{teo}
The paper is organized as follows. In section 2, we fix notation and
recall some general facts. In Section 3, we show the main tools that
we will use to determine the group $\Aut(X_0^*(N))$ for a {{fixed}}
value of $N$. Sections 4 and 5 are devoted to prove Theorem
\ref{main2} for odd and even levels respectively. The key point in
these two last sections is to delimit a finite  set of positive
integers containing all levels $N$ such that $X_0^*(N)$ has a non
trivial automorphism group.  In section 4, for the odd levels, this
result follows from Proposition 3. This proposition is  based on  an
idea already used in \cite{KeMo}[Lemma 2.7], \cite{BH}[Lemma 3.3]
and \cite{Go16}[Lemma 6], but always proved with different arguments
because the modular curves involved  in these statements are
different. In section 5, the delimitation of a finite set for the
even levels  is obtained from  Proposition 4. This proposition
presents   an unknown inequality involving the genera of the curves
$X_0^*(N)$ and $X_0^*(2N)$ for all odd square-free values of  $N$.
\section{Notation and general facts}

Let $N>1$ be a square-free integer. We fix, once and for all, the following
notation.
\begin{itemize}
\item [(i)] We denote by  $B(N)$  the group of the Atkin-Lehner involutions of $X_0(N)$. If $N'|N$, $B(N')$  can be also considered as  the subgroup of $B(N)$ formed by the Atkin-Lehner involutions $w_d$ such that $d|N'$.
    \item [(ii)] The integer  $\omega(N)$ is the number of primes dividing $N$. In particular, $|B(N)|=2^{\omega(N)}$.
\item [(iii)] The integers
 $g_N$ and $g_N^*$ are the genus of $X_0(N)$ and $X_0^*(N)$ respectively.

\item[(iv)]
 We denote by $\New_N$ the  set of normalized newforms in $S_2(\Gamma_0(N))$ and $\New_N^*$ is the subset of $\New_N$ formed by the newforms  invariant under the action of the group  of the Atkin-Lehner involutions $B(N)$.

\item[ (v)]  $S_2(N)=S_2(\Gamma_0(N))$,  $S_2^*(N)$ is the vector space $S_2(N)^{B(N)}$, $J_0(N)=\Jac(X_0(N))$ and $J_0^*(N)=\Jac(X_0^*(N))$.

\item [(vi)] For a normalized  eigenform $h\in S_2(\Gamma_0(N))$, $A_h$ denotes the abelian variety defined over $\Q$ attached by Shimura to $h$. This abelian variety can be viewed as the optimal quotient of $J_0(N)$ such that the pullback of $\Omega^1_{A_h}$ is the vector space generated by the Galois conjugates of $h(q)\,dq/q$.

\item[ (vii)] Given two abelian varieties $A$ and $B$ defined over the number field $K$, the notation $A\stackrel{K}\sim B$ stands for $A$ and $B$ are isogenous over $K$.

\item [(viii)] For an integer  $m\geq 1$ and $f\in\New_N$, $a_m(f)$  is the $m$-th Fourier coefficient of $f$.

\item[(ix)] As usual, $\psi$ denotes the Dedekind psi function.

\item[(x)] $G_\Q$ denotes the absolute galois group $\Gal (\overline{\Q}/\Q)$ once an algebraic closure of $\Q$ has been fixed.
\end{itemize}

 We recall that $\Omega^1_{X_0(N)/\Q}$  is the $\Q$-vector subspace of elements in  $S_2(N)dq/q$ with rational $q$-expansion, i.e. $S_2(N) dq/q\cap\Q[[q]]$.  For $N$ is square-free and $g_N^*\geq 1$, it is known that
$$
J_0^*(N)\stackrel{\Q}\sim \prod_{1<M|N}\prod _{f\in \New_M^*/G_{\Q}}A_f\,.
$$
These abelian varieties $A_f$  in the decomposition are simple and pairwise non-isogenous over $\overline{\Q}$ and the endomorphism algebra $\End(\Jac(X_0^*(N))\otimes \Q$ is isomorphic to the product of totally real numbers fields (cf. \cite[\S 2]{BH}). A basis of $S_2^*(N)$ is formed by the eigenforms
$$\bigcup_{1<M|N} \{ \sum_{1\leq d|N/M} d\, f(q^d)\colon f \in\New_M^*\}\,.$$
Every  $f_i\in\New_{M_i}^*$ with
 $M_i|N$ determines the normalized eigenform $h_i=\sum_{1\leq d|N/M_i} d\, f(q^d)$ in $S_2^*(N)$ such that $A_{f_i }\stackrel{\Q}\sim  A_{h_i}$.
 On the one hand, the basis of the Galois conjugated of the newforms $f_i$   allows us to compute $|X_0^*(\F_p^n)|$ for all primes $p\nmid N$,  thanks to the Eichler-Shimura congruence. Indeed, the characteristic polynomial of $\Frob_p$ acting on the Tate module of $J_0^*(N)$ is
 $$  P(x)=\prod_{1>M|N}\prod_{f\in \New_M^*}(x^2-a_f(f)x+p)=\prod_{i=1}^{2\,g_N^*}(x-\alpha_i)\,,$$
 and, thus,
 $$|X_0^*(\F_{p^n})|=p^n+1-\sum _{i=1}^{2\,g_N^*}\alpha_i^n\,.$$
 On the other hand, the basis of the  regular differentials formed by all Galois conjugates of $h_i(q) \, dq/q$ allows us  to compute equations for $X_0^*(N)$ when $g_N^*>1$.

 Similarly, if $X$ is a curve defined over $\Q$  of genus $g>0$ for which there is a non constant morphism $X_0^*(N)\rightarrow X$, then there is a subset  $\cS$ of $\cup_{1<M|N}\New_M^*$ stable by Galois conjugation such that
 $$
 \Jac (X)\stackrel{\Q}\sim \prod_{f\in\cS/G_\Q}A_f\,
 $$
Hence, for a prime $p\nmid N$, the characteristic polynomial of $\Frob_p$ acting on the Tate module of $\Jac (X)$ is
 $$  P(x)=\prod_{f\in \cS}(x^2-a_f(f)x+p)=\prod_{i=1}^{2\,g}(x-\alpha_i)$$
 and  $|X(\F_{p^n})|=p^n+1-\sum _{i=1}^{2\,g}\alpha_i^n$.  The set of the normalized eigenforms atached to $\cS$ in $S_2^*(N)$ is a basis of $\pi^*(\Omega^1_X) q/dq$ that allows us  to compute equations for $X$ when $g>1$.

\section{Preliminary results}
In this section, we present some tools that will be applied to determine  wether  the group  $\Aut(X_0^*(N))$ is trivial or not, and to determine exactly this group when it is non trivial.

For a curve $X$ defined ove $\Q$ of genus $\geq 2$ and of good reduction at a prime $p$,  one has that $\Aut_{\Q}(X)\hookrightarrow \Aut_{\F_p}(X\otimes\F_p)$. We can apply Theorem 2.1 in  \cite{Go17} to discard the existence of involutions defined over $\Q$. More precisely, we will use the following criterion.
\begin{lema}\label{criterion1}
Let $X$ be a curve defined over $\F_p$ of genus $g>2$. Consider
the sequence
$$P_p(n):= \operatorname{mod\,} [(\sum_{d|n}\mu(n/d)|X(\F_{p^d})|)/n,2]\, $$
 where  $\operatorname{mod\,} [r,2]$ denotes $0$ or $1$ depending on whether $r$ is even or not, and
$\mu$ is the {Moebius} function.
If there is an integer $k>0$ such that
 $\sum_{n\geq 0}^k(2 n+1) P_p(2n+1)> 2 g+2$, then  $\Aut_{\F_p}(X)$ does not contain any involution.
\end{lema}
We recall that, for a nonhyperelliptic curve $X$ defined over $\C$
with genus $g>2$, the image of the canonical map $ X\rightarrow
\PP^{g-1}$ is the common zero locus of a set of homogeneous
polynomials of degree $2$  and $3$, when $g>3$,  or of a homogenous
polynomial of degree $4$, if $g=3$.

More precisely, assume that  $X$ is defined over $\Q$ and  choose a basis $\omega_1,\cdots,\omega_g$ of $\Omega^1_{X/\Q}$.
For any integer $i\geq 2$, let us denote by $\cL_i$ the $\Q$-vector space of homogeneous polynomials $Q\in\Q [x_1,\cdots,x_g]$
of  degree $i$ that satisfy $Q(\omega_1,\cdots,\omega_g)=0$. Of course,  $\dim \cL_i\leq \dim \cL_{i+1}$ because one has $x_j\cdot Q \in\cL_{i+1}$ for all $Q\in
\cL_i$  and for $1\leq j\leq g$.

If $g=3$, then $\dim \cL_2=\dim \cL_3=0$ and $\dim \cL_4=1$. Any generator of $\cL_4$ provides an equation for $X$.
For $g>3$, $\dim \cL_2=(g-2)(g-3)/2>0$ and a basis of $\cL_2\bigoplus\cL_3$ provides a system of equations for $X$.
When $X$ is neither trigonal nor a smooth plane quintic ($g=6$), it suffices to take a basis of  $\cL_2$.

As said in the above section,
$J_0^*(N)\stackrel{\Q}\sim A_{h_1}\times \cdots\times
A_{h_n}$ for some normalized
eigenforms $h_1,\cdots,h_k\in S_2^*(N)$.  These abelian
varieties are simple and pairwise nonisogenous over $\Q$ and, any
involution $u$ of the curve leaves stable each $A_{h_i}$  acting on
$\Omega^1_{A_{g_i}}$ as the product by $-1$ or the identity.

Choose a basis $\{\omega_1,\cdots,\omega_{g_N^*}\}$ of $\Omega^1_{X_0^*(N)/\Q}$ obtained as  the ordered  union
of bases of all $\Omega^1_{A_{h_i}/\Q}$.  An involution $u$ of $X_0^*(N)$ induces a linear map $u^*:\Omega^1_{X_0^*(N)/\Q}\rightarrow \Omega^1_{X_0^*(N)/\Q}$ sending $(\omega_1,\cdots,\omega_{g_N^*})$ to $(\varepsilon_1\omega_1,\cdots,\varepsilon_n\omega_{g_N^*})$ with $\varepsilon_i=\pm 1$  for all $i\leq {g_N^*}$ and satisfying
\begin{equation}\label{involution}
Q(\varepsilon_1x_1,\cdots,\varepsilon_{g_N^*}  x_{g_N^*})\in \cL_i \text{ for all } Q\in\cL_i \text{ and for all } i\geq 2\,.
\end{equation}
 The genus of the quotient curve is the cardinality of the set $\cI= \{i\colon \varepsilon_i=1\}$ and $\{\omega_j\}_{ j\in\cI}$ is a basis of the pullback of the regular differentials of the quotient curve.
Conversely, for a linear map $u^*$ as above satisfying condition (\ref{involution}), only one of the two maps $\pm u^*$ comes from an involution of the curve, because we are assuming that $X$ is nonhyperelliptic.

 \vskip 0.2 cm
 We particularize this fact to our case, that will be the main tool to determine the group $\Aut(X_0^*(N))$ for a fixed level $N$.
\begin{lema}\label{con}
 Assume $X_0^*(N)$ is nonhyperelliptic, i.e. $g_N^*>2$. Let $\omega_1,\cdots,\omega_{g_N^*}$ be a basis of $\Omega^1_{X_0^*(N)/\Q}$ as above. Then,
  \begin{itemize}
  \item [\rm{(i)}]If $u$ is a non trivial involution of $X_0^*(N)$ such that the  genus of the curve $X_u=X_0^*(N)/\langle u\rangle$ is $g_u$  and $\Jac (X_u)\stackrel{\Q}\sim A_{h_1}\times \cdots \times A_{h_k}$ ($k<n$), then
\begin{equation} \label{eq}Q(-x_1,-x_2,\cdots,-x_{g_u}, x_{g_u+1},\cdots x_{g_N^*-1}, x_{g_N^*})\in \cL_i \text{ for all } Q\in\cL_i \text{ and for all } i\geq 2\,.\end{equation}
\item [\rm{(ii)}]If $\dim A_{h_1}\times \cdots \times A_{h_k}=g_u$ and the condition (\ref{eq}) is satisfied, then there exists an involution  $v$ of $X_0^*(N)$ such that $$\Jac (X/\langle v\rangle )\stackrel{\Q}\sim A_{h_1}\times \cdots \times A_{h_k} \text{ or }\Jac (X/\langle v\rangle)\stackrel{\Q}\sim A_{h_{k+1}}\times \cdots\times A_{h_n}\,$$
\end{itemize}
\end{lema}
In order to make easy the computation of the condition (\ref{eq}) for $\cL_2$, we introduce the vector subspace $\cL_2^{ns}$ of polynomials in $\cL_2$ that do not contain square monomials, i.e. polynomials of the form $\sum_{1\leq i <j\leq g}a_{ij}x_i x_j$. For any $Q\in \cL_2$ and  any $r\leq g$, we have that
\begin{equation}\label{eq2}
Q(-x_1,\cdots,-x_r,x_{r+1},\cdots,x_g)\in\cL_2\Leftrightarrow Q-Q(-x_1,\cdots,-x_r,x_{r+1},\cdots,x_g)\in\cL_2^{ns}\,.
\end{equation}
In general, it is expected  that $\dim \cL_2^{ns}=\operatorname{Max}(\dim \cL_2-g,0)=\operatorname{Max}((g-1)(g-6)/2,0)$. When $\dim \cL_2^{ns}=0$, the condition (\ref{eq2}) amounts to saying that $Q$ is simultaneously even in the variables $x_1,\cdots,x_r$.

\begin{rem}\label{nopetri}
For a high genera $g_N^*$, if we can not to claim that the automorphism group of  $X_0^*(N)$ is   trivial by means of Lemma \ref{criterion1},
we have some aternative ways to prove this fact,  without applying   Lemma \ref{con}  to this curve. Indeed, let
 $\omega_1,\cdots,\omega_{g_u}\in \Omega^1_{X_0^*(N)/\Q}$ be a  basis of the pullback of regular differentials  of a possible quotient $X_u$. If for a prime $\ell\nmid N$, $|X_0^*(N)(\F_{\ell^n})|> 2 |X_u(\F_{\ell^n})|$, we can discard this  possibility. Another way is the following. If $g_u>2$ and we know
 that $X_u$ is not  hiperelliptic, we apply  Petri Theorem to regular differentials $\omega_1,\cdots,\omega_{g_u}$. If $g_u>3$ and the  corresponding vector space    $\cL_2$  does not have dimension $(g_u-3)(g_u-2)/2$, then  we can also  discard this  possibility. For  $g_u=3$ and $\dim \cL_4\neq 1$, we can proceed in the same way.

\end{rem}

\vskip 0.2 cm
The following result allows us to use  the knowledge of the group of automorphisms  of lower levels,  that can be a source for getting bounds for the even levels $N$ such  that the curves  $X_0^*(N)$ may have non trivial automorphisms.

\begin{prop}\label{restrict}
Let $N>1$ be a square-free integer and let  $p$ be a prime dividing $N$ such that  $g_{N/p}^*> 1$. It there exists a non trivial involution $u$ of $X_0^*( N)$ such that  the action of $u$ on $J_0^*(N/p)$ is non trivial, then the group $\Aut(X_0^*(N/p))$ is non trivial or the curve  $X_0^*(N/p)/ \F_p$ is hyperelliptic.

\end{prop}
\noindent
{\bf Proof.}  Let $\phi$ and $\phi_p$ be the morphisms from $X_0(N)$ to $X_0(N/p)$ induced by the automorphisms of the complex upper half-plane given by $z\mapsto z$ and $z\mapsto p z$ respectively. Consider the morphism $\nu=\phi^*+p\,
\phi_p^* \colon J_0(N/p)\rightarrow J_0( N)$ that is defined over $\Q$. Since for every cusp
form $h\in S_2^*(N/p)$ one has that $\nu (h)\in S_2^*(N)$,  then $\nu$ induces a morphism from $J_0^*(N/p)$ to $J_0^*(N)$ that will  still be denoted by $\nu$.
It is clear that $J_0^*(N/p)$  and $\nu(J_0^*(N/p))$ have the same dimension. Moreover, for any Hecke operator $T_M$ with $\gcd(M,N/p)=1
$, one has $\nu\cdot T_M=T_M\cdot \nu$ (cf. \cite{Li}). Recall that $\End J_0^*(N/p)$ and $\End J_0^*(N)$ are generated by Hecke operators (cf. \cite{Ribet}[Proposition3.2]).

  Let $u$ be of a non trivial involution of $X_0^*( N)$. This involution induce an involution  $u^*$ of $J_0^*( N)$ which left stable $\nu(J_0^*(N/p))$ because $\nu(J_0^*(N/p))$ and $J_0^*(N)/\nu(J_0^*(N/p))$ are defined over $\Q$ and do not have isogenous quotients. Hence $u^*$ can be viewed as an  involution of $J_0^*(N/p)$.  Now, assume that this one is non trivial. On the one hand, the normalization of $X_0^*(N)/ \F_p$ is  the curve $X_0^*(N/p)/\F_p$ (cf. \cite{Ha97}[Section 5]) and $u$ induces an involution $\widetilde u$ of $X_0^*(N/p)/\F_p$. On the other hand,
    the reduction of $u^*$  modulo $p$ is a  non trivial  involution of $J_0^*(N/p)/\F_p$ that  coincides with the endomorphism induced by the involution $\widetilde u$.  We claim that if $\widetilde u $ is not hyperelliptic, then $\widetilde u$ raises to an involution of $X_0^*(N/p)$.
    Indeed,  $\widetilde u$  respects the canonical polarization of the curve  $X_0^*(N/p)/\F_p$. Since the reduction of the endomorphism ring is injective, this equality is still true in $ \Jac (X_0^*(N/p))/\Q $ and, thus, $u^*$  is induced by an automorphism of the curve $X_0^*(N/p)/\Q$ by Torelli's Theorem.
   \hfill $\Box$
\begin{rem}\label{restgenus}
 Assume that  there exist a non trivial involution $u$ of $X_0^*( N)$  and a prime $p|N$ with $g_{N/p}^*>2$.  If  the group $\Aut(X_0^*(N/p))$ is trivial and $X_0^*(N/p) / \F_p$ is not hyperelliptic, then $g_u \geq g_{N/p}^*$ and the pullback of  $\Omega^1_{J_0^*(N/p)}$ in $\Omega^1_{J_0^*(N)}$ (under $\nu$) is contained in the pullback of $\Omega^1_{\Jac (X_u)}$. In particular,  $ g_{N/p}^*\leq (g_N^*+1)/2$ because  by Hurwitz one has $ g_u\leq (g_N^*+1)/2$.
\end{rem}

Note that a necessary condition to be  $X_0^*(N/p)/\F_p$  hyperelliptic  is that the following inequa\-li\-ties are satisfied
 $$|X_0^*(N/p)(\F_{p^n})|\leq 2 p^n+2 \text{ for all integers $n>0$.} $$
 The hyperelliptic curves have a different behaviour depending on wether are defined over fields od characterisitc $2$ or not.
 Since  \cite{BGGP}[Lemma 2.5] can be applied to hyperelliptic curves defined over fields of characteristic different  from $ 2$,  we have a similar result to  \cite{BGGP}[Lemma 6] for modular curves over $\F_p$ with $p\neq 2$. More precisely, in our case:
\begin{lema}\label{hip}
   Assume $p$ is an odd prime not dividing  $M$.  Let $X/\Q$ be a curve of genus $g>2$ for which there is a no constant morphism $\pi:X_0^*(M)\rightarrow X$ and let $\cS$ the $\Q$-vector subspace  $\pi^*(\Omega^1_{X/\Q} ) q/dq$ of $S_2^*(M)$. The curve   $X/ \F_p$ is hyperelliptic if, and only if, there is a basis $f_1,\cdots, f_{g}$ of $\cS$ with integer $q$-expansions whose reductions mod $p$ satisfy
\begin{equation}\label{congq}
f_i(q) \,(\operatorname{mod}\, p)=\left\{\begin{array}{cr} q^i+ O(q^i)&\text{ if the cusp $\infty$ is not a Weierstrass point of $X/\F_p$,}\\q^{2 i-1}+ O(q^{2i-1})& \text{otherwise.}\end{array} \right. \,,
\end{equation}
and  for any  such a basis,  the functions on $X/\F_p$ defined by
$$
x=\frac{f_{g-1}}{f_g}\pmod p\,,\quad  y=\frac{ q dx/dq}{f_{g}}\pmod p\,,
$$
satisfy $y^2=P(x)$ for  an unique square-free polynomial $P(X)\in\F_p[X]$ which has  degree $2g+1$ or $2g+2$ depending on whether the cusp $\infty$ is a Weierstrass point or not.

\end{lema}

For the case $p=2$, we  use the following result.

\begin{lema}\label{hyp2} Let $N>1$  be an odd square-free integer such that $g_N^*>2$. If $X_0^*(N)/\F_2$ is hyperelliptic, then $N$ is in the set $\{183,185,187,203,335,345,385\}$.

\end{lema}
\noindent{\bf Proof.} Put $n=\omega(N)$.  Again by Ogg, if $X_0^*(N)/\F_2$ is hyperelliptic, then
$$
\frac{\psi(N)}{12}+ 2^{n}\leq 2^{n+1}|\PP^1(\F_4) |=5\cdot 2^{n+1}\,.
$$
Hence, $\psi (N)\leq 2^ n 108$ and, after discarding the values $N$ such that $g_N^*\leq 2$, we obtain the following $59$ values of $N$:
$$
\begin{array}{r}
97, 109, 113, 127, 137, 139, 149, 151, 157, 163, 173, 179, 181, 183,
185, 187, 193, 197, 199, 201, 203, \\
211, 217, 219, 235, 237, 247, 249,
253, 259, 265, 267, 273, 291, 295, 301, 303, 305, 309, 319, 321, 323, \\
329, 335, 341, 345, 355, 371, 377, 385, 391, 399, 429, 435, 455, 465,
483, 561, 595\,.
\end{array}
$$
We can discard the levels $N$ such that $A_{N,m}:=|X_0^*(N)(\F_{2^m})| - 2(2^m+1)>0$. Since
$$
\begin{array}{|c|c|}
(N,m) &  A_{N ,m}\\\hline
(97,1) &1\\\hline
(109,2) & 1\\\hline
(113,2) &1\\ \hline
(137,2) &2\\ \hline
(139,2) &1\\ \hline
(149,2)&2\\ \hline
(151,2)&1\\ \hline
(157,1) &2\\ \hline
(163,1)&2\\ \hline
(173,2 &4\\ \hline
(179,2) &2\\\hline
(181,2)&4\\ \hline
\end{array}\quad
\begin{array}{|c|c|}
(N,m) &  A_{N ,m}\\\hline
(193,1)& 2\\ \hline
(197,2)& 5\\ \hline
(199,2)&2\\\hline
(201,1)&1\\ \hline
(211,2)&5\\ \hline
(219,1)&1\\ \hline
(235,1)&1\\ \hline
(237,1)&1\\ \hline
(249,2)& 4\\ \hline
(253,1)&1\\ \hline
(259,2)&2\\ \hline
(265,1)&2\\ \hline
\end{array}
\quad
\begin{array}{|c|c|}
(N,m) &  A_{N ,m}\\\hline
(267,2)& 4\\ \hline
(273,1)&1\\ \hline
(291,1)& 3\\ \hline
(295,2)&2\\ \hline
(301,2)&3\\ \hline
(303,2)&3\\ \hline
(305,2)&4\\ \hline
(309,1)&2\\ \hline
(319,2)&1\\ \hline
(321,2)&5\\ \hline
(323,2)&2\\ \hline
(329,2)& 2\\ \hline
\end{array}
\quad
\begin{array}{|c|c|}
(N,m) &  A_{N ,m}\\\hline
(341,2)&5\\ \hline
(355,2)&3\\ \hline
(371,2)& 5\\ \hline
(377,2)&4\\ \hline
(391,1)& 1\\ \hline
(399,1)&1\\ \hline
(429,2) &1\\ \hline
(435,1)& 1\\ \hline
(465,2)& 2\\ \hline
(483,2)& 5\\ \hline
(561,1)&1\\ \hline
(592,2)& 3\\ \hline
\end{array}
$$
By Lemma \ref{criterion1}, we also can discard the values $N$ for which we know that $X_0^*(N)/\F_2$ does not have any involution:
$
127,217,247
$. \hfill $\Box$

Note that for a prime $p$ we know that the automorphism group is
trivial when $g_p^*>2$ (cf. \cite{BH}[Theorem 1.1]). For  $g_N^*=3$,
the group $\Aut (X_0^*(N))$ is non trivial if, and only if,
$X_0^*(N)$ is bielliptic (cf. \cite{BaGon}{Lemma 13]). For this
reason in the next sections, we exclude the cases $g_N^*\leq 3$ that
can be found in Table 4 in \cite{BaGon2}[ Appendix].

\section{Odd levels}

We know  that when $g_N^*>1$, if   there is a non trivial involution $u$ of $X_0^*(N)$, then the cusp $\infty$ is not a fixed point of $u$ (cf. \cite{BH}[Lemma 3.2]). The following result, similar to Lemma 6 in \cite{Go16}, is the key result that allow us  to delimit a finite set of  odd square-free integers  containing all odd square-integers $N$ such that the group $\Aut (X_0^*(N))$ is non trivial.

\begin{prop}
Assume that  the square-free $N$ is odd, $g_N^*>1$ and there is a non trivial involution $u$ of $X_0^*(N)$. Then, the $\Q$-gonality of $X_0^*(N)$ is $\leq 6$ and $u$ has at most $12$ fixed points.
\end{prop}
\noindent
{\bf Proof.}  Put $Q=u(\infty)$ and let $P\in X_0(N)$ be such that $\pi (P)=u(\infty)$, where $\pi $ is the natural projection of $X_0(N)$ onto $X_0^*(N)$.
Since $Q$ is not a cusp, there is
an elliptic curve $E$  defined over $\overline {\Q}$ and a $N$-cyclic subgroup $C_N$ of $E(\overline {\Q})$ such that $P=(E,C_N)$. The other pre-images of $Q$ under $\pi$ are the points
$$w_{d}(P)=(E/C_d, (E[d]+ C_{N/d})/C_d)\quad \forall d|N\,,$$
where $C_{d}$ denotes the $d$-cyclic subgroup of $C$, that is $C_{d}=C_N\cap E[d]$.

 For any noncuspidal $S\in X_0^*(N)$, we consider the divisor
$$D_S:=(u T_2 -T_2 u)\left(\infty-S\right)\,,$$
where $T_2$ denotes the Hecke operator  viewed as a  correspondence of the curve $X_0^*(N)$. We claim that $D_S$ is nonzero but linearly
equivalent to zero.

The endomorphism algebra  $\End (J_0^*(N))\otimes \Q$ is  commutative. Hence,   $u\, T_2=T_2\, u$
and  $D_S$ is a  principal divisor.

 Next, we will prove  that $D_S$ is nonzero.
 If $D_S$ is a  zero divisor, then  $u T_2 (\infty)$ must be equal to $T_2 u(\infty)$ because $T_2(\infty)=3\infty$ and $\infty$ is not in the support of $T_2(S)$.  To prove that $D_S$ is a nonzero divisor,  we only need to prove that the condition $3(Q)=T_2(Q)$ cannot occur for a noncuspidal point  $Q\in X_0^*(N)$.

Let $G_i$, $1\leq i\leq 3$, be the three $2$-cyclic subgroups of $E [2]$. Since
$$T_2(Q)=\sum_ {i=1}^3 \pi ((E/G_i, (C_N +G_i)/G_i))\,,
$$ the condition $3(Q)= T_2(Q)$  implies that each elliptic curve $E/G_i$ is isomorphic to $E/C_d$ for some  $d|N$. Therefore, $E$ has an endomorphism whose kernel is a $2d$-cyclic subgroup and, thus,  $E$ has CM by a quadratic order $\cO$ of discriminant $D$.
The conductor of the  discriminant $D$ cannot be even because  $2d$ is a norm of $\cO$ and $2d\not\equiv 0\pmod 4$.
  Since $\pi (P)=\pi (w_d(P))$, this property holds for all elliptic curves $E/C_d$ and, thus, also for all elliptic curves $E/G_i$.

Now, we claim that for every elliptic curve $E$ with CM by the order of discriminant  $D$  with odd conductor,  there is at least a $2$-subgroup $G$ of $E[2]$ such that the discriminant of the order $\End(E/G)$  has even conductor. This fact concludes that $T_2(Q)\neq 3Q$ and, thus, $D_S$ is nonzero. Indeed, let  $[a,b,c]:=ax^2+bxy+c y^2\in \Z[x,y]$ be a primitive quadratic form of discriminant $D$ (with odd conductor). The primitive quadratic forms attached to the three elliptic curves $E/G_i$  are
$$Q_1=[4a,2b,c]/\gcd(c,2),\,\, Q_2=[a,2b, 4c]/\gcd(a,2),\,\, Q_3=[4a, 2b-4a,a-b+c]/\gcd(2,a-b+c)\,.
$$
  If the discriminants
  of $Q_1$ and $Q_2$ are equal to $D$, then $a$ and $c$ must be  even.   Since  $[a,b,c]$ is primitive, $b$ must be  odd and this fact leads to the contradiction that the discriminant  of $Q_3$ is $4D$, with even conductor.

By taking $S=u(\infty)$, $D_S$ is defined over $\Q$ and, thus, the $\Q$-gonality is at most $6$. Finally,  since $u^*(D_S)\neq D_S$ for some noncuspidal point $S\in X_0^*(N)(\C)$,   any  non trivial automorphism of   $X_0^*(N)$ has at most $12$ fixed points (cf. Lemma 3.5 of \cite{BH}).
\hfill $\Box$

\begin{cor}\label{trivial}
When $N$  is odd, if $g_u$ is the genus of the curve  $X_0^*(N)/\langle u\rangle$ for a non trivial involution $u$ of $X_0^*(N)$, then
$$  \frac{g_N^*-5}{2}\leq  g_u\leq \frac{g_N^*+1}{2}\,.$$
Moreover, if $J_0^*(N)$ has a simple factor of dimension larger than $(g_N^*+5)/2$, then $\Aut(X_0^*(N))$ is trivial.
\end{cor}

\noindent{\bf Proof.} The inequalities follows from the Riemann-Hurwitz formula applied to the projection $X_0^*(N)\rightarrow X_0^*(N)/\langle u\rangle$. For the last assertion see \cite{BH}[Corollary 3.9] \hfill $\Box$

\begin{lema} Assume that  the $\Q$-gonality of $X_0^*(N)$ is at most $6$.  If $N$ is odd, then
\begin{equation}\label{gonality2}\psi(N)\leq 2^{\omega(N)} 348\,.\end{equation}
\end{lema}
\noindent{\bf Proof.} We know that $\displaystyle{ 2^{\omega(N)}+\frac{\psi(N)}{12}\leq |X_0(N)(\F_4)|}$ (cf. proof  of \cite{Ogg}{Theorem 3]). The statement follows from the fact that $|X_0(N)(\F_4)|\leq 2^{\omega(N)}\cdot6\cdot |\PP^1(\F_4)|$.\hfill $\Box$

There are $471$ values of $N$  (square-free and odd, {{$N\geq3$}})
satisfying the condition (\ref{gonality2}), whose ma\-ximum value is
$3003$. Excluding all values that are prime or $g_N^*\leq 3$, we
obtain $293$ values.
  By applying Lemma \ref{criterion1} to $X_0^*(N)/\F_p$, we can discard  $248$ values of $N$. More precisely, with $p=2$ the values
 $$ \begin{array}{l}247, 253, 259, 267, 291, 301, 305, 319, 323, 327, 339, 355, 365, 377, 381, 391, 393, 395, 403, 407, 411, \\
 413, 417, 427, 451, 453, 469, 471, 473, 481, 485, 489, 493, 501, 505, 511, 519, 527, 533, 535, 537, 543,\\
  545, 553, 559, 565,  573, 579, 581, 589, 591, 595, 597, 611, 627, 629, 633, 635, 651, 655, 667, 669, 671, \\
  679, 681, 685, 687, 695, 697, 699, 703, 707, 713, 717, 721, 723, 731, 737, 741, 745, 749, 755, 759, 763,\\
771, 777, 779, 781, 785, 789, 791, 793, 795, 799, 803, 805, 807, 813, 815, 817, 831, 835, 843, 849, 851,\\ 865,
 869, 871, 879, 885, 889, 893, 895, 897, 899,901, 903, 905, 913, 915, 917, 921, 923, 933, 935, 939,\\ 943, 949,
  951, 955, 959, 965, 969, 973, 979, 985, 993, 995, 1001, 1003, 1005, 1007, 1011, 1015, 1023,\\ 1027, 1037, 1041,
  1043, 1045, 1057, 1065, 1067, 1073, 1081, 1085, 1095, 1099, 1105, 1111, 1113, 1115,\\1121, 1131, 1133, 1135,
  1139, 1141, 1145, 1147, 1157, 1159, 1169, 1173, 1177, 1185, 1189, 1199, 1207, \\1209, 1211, 1219, 1221, 1235,
  1239, 1241, 1243, 1245, 1247, 1261, 1265, 1271, 1273, 1281, 1295, 1311, \\1353, 1407, 1419, 1435, 1443, 1455,
  1463, 1479, 1491, 1505, 1515, 1533, 1545, 1547, 1581, 1595, 1599, \\1605, 1635, 1645, 1653, 1659, 1677, 1695,
  1705, 1729, 1743, 1749, 1767, 1771, 1785, 1833, 1855,1885, \\1887, 1955, 1995, 2015, 2035, 2093, 2145, 2415, 2805, 3003\,,
\end{array}
$$
with $p=3$ the values $445,1495,1615$, with $p=5$ the values $623$,
with $p=7$ the values $583,753,1551$  and  for $p=11$ the value
$1335$. {{Corollary \ref{trivial} allows us to exclude the values
  $$
     235, 237, 273,341,385, 415,435, 497,515,517,649,767,715, 989,1079,1309 \,.
 $$}}

  So, we only have to consider  $29$ values for $N$, that we present  together the splitting of the corresponding jacobians collected by the genus $g_N^*$.
  Fron know on, the splitting $$J_0^*(N)\stackrel{\Q}\sim \prod_{i=1}^r A_{f_i}\,, \text{  with } f_i\in \New_{M_i}^* \text{  and }\dim A_{fi}=n_i$$ will be presented as ${n_1}_{M_1}+\cdots+{n_r}_{M_r}$. Obviusly,  $g_N^*=\sum_{i=1}^r n_i$ and, for a divisor $M$ of $N$, one has $g_M^*=\sum_{M_i|N}n_i$.

  $$
\begin{array}{|c|r||}
 N& J_0^*(N)\\ \hline\hline
201& 2_{67}+1_{201}+1_{201} \\[ 3 pt]
219 & 2_{73}+1_{219}+1_{219}   \\[ 3 pt]
321&2_{107}+2_{321}\\[ 3 pt]
335&2_{67}+2_{335}\\[ 3 pt]
345&2_{115}+2_{345}\\[ 3 pt]
399& 1_{57}+2_{133}+1_{399}\\[2 pt]
483 &2_{161}+2_{483}
\\ \hline
\hline
371&1_{53}+1_{371}+3_{371} \\[ 3 pt]
465& 2_{93}+1_{155}+2_{465}   \\[ 3 pt]
551&2_{551}+3_{551}\\[ 4 pt]
555&1_{37}+1_{185}+1_{185}+2_{555}\\[ 4 pt]
645&1_{43}+1_{129}+1_{215}+2_{645}\\[ 4 pt]
663& 2_{221}+3_{663}\\[4 pt] &  \\\hline
\end{array}
\begin{array}{|c|r||}
 N& J_0^*(N)\\ \hline\hline
 265&1_{53}+1_{265}+2_{265}+2_{265} \\
447& 3_{149}+3_{447}   \\
561&3_{187}+3_{561}\\
609&3_{203}+3_{609}\\
615&1_{123}+2_{205}+1_{615}+2_{615}\\ \hline\hline
 309&2_{103}+5_{309} \\
437& 2_{437}+5_{437}   \\
861&1_{123}+2_{287}+4_{861}\\
\hline
\hline
665&2_{133}+6_{665}\\
1155&1_{77}+2_{165}+3_{385}+1_{1155}+1_{1155}\\
\hline\hline
689& 1_{53}+1_{689}+2_{689} +2_{689}+3_{689}\\
705 & 1_{141}+5_{235}+1_{705}  +2_{705}\\
987&1_{141}+3_{329}+2_{987}+3_{987}\\
1365&1_{65}+1_{91}+3_{273}+1_{455}+3_{1365}\\ \hline\hline
957& 4_{319}+7_{957}\\ \hline\hline
1055 & 3_{211}+3_{211}+3_{1055}  +6_{1055}\\\hline
 \end{array}
  $$

\vskip 0.2 cm

\begin{prop}\label{odd}
Let $N>1$ be an odd square-free  integer such that $g_N^*>3$. The group $\Aut (X_0^*(N))$ is non trivial if, and only if, $N=645$. In this case, the order of this group is $2$ and an equation  for the quotient curve by the non trivial involution is given by
$$Y^2=X^6 + 8 X^4 + 20 X^2 + 12 X + 4\,.$$
\end{prop}
\noindent {\bf Proof.}
To apply Petri's Theorem to the curve $X_0^*(N)$, we will  take  a basis of $\Omega^1_{X_0^*(N)/\Q}$ as in Lemma \ref{con}, following the order showed in the splitting tables.

 For $g_N^*=4$, $\dim \cL_2=1$. In all cases in the above table,  $\dim \cL_2^{ns}=0$. The genus of a quotient curve by an involution  must be $2$. We have to check the condition (\ref{eq}) for all pairs $x_i,$ and $x_j$ such that $\omega_i$ and $\omega_j$ are a basis of the regular differentials  corresponding to an abelian quotient of $J_0^*(N)$ of dimension $2$.
 A nonzero polynomial $Q\in\cL_2$ neither  satisfies $Q(-x_1,-x_2,x_3,x_4)=Q$ for $N\neq 399$ nor $Q(x_1,-x_2,-x_3,x_4)=Q$ for $N=399$. Therefore, for all  these cases the curves $X_0^*(N)$ have trivial automorphism group.

For  $g_N^*=5$, $\dim \cL_2=3$. The genus of a quotient curve by an involution  must be $2$  or $3$. It is suffices to check the condition (\ref{eq}) for all pairs $x_i,$ and $x_j$ such that $\omega_i$ and $\omega_j$ is a basis of the regular differentials  corresponding to an abelian quotient of $J_0^*(N)$ of dimension $2$.
 For the values in the above table $N\neq 645$, one has $\dim \cL_2^{ns}=0$ and the polynomials of $\cL_2$ are not simultaneously even in the variables $x_i,x_j$. Hence, for all  these cases the curves $X_0^*(N)$ have trivial automorphism group. For $N=645$, we know that $X_0^*(645)$ is not trigonal (cf. \cite{HaSh}) and, thus,  $\cL_2$ defines the curve.  The set $\cup_{M|645}\New_M^*$  is
 $$ \New_{43}^*=\{f_1\},\quad \New_{129}^*=\{f_2\}, \quad\New_{215}^*=\{f_3\},$$
 $$ \New_{645 }^*=\{f_4=q+\sqrt 2 q^2+\cdots,f_5=q-\sqrt 2 q^2+\cdots\}\,.$$
 Taking the following basis of $\Omega^1_{X_0^*(645)/\Q}$:
 $$
 \omega_1=(\sum_{d|15} d \,f_1(q^d))dq/q\,, \,\omega_2=(\sum_{d|5} d \,f_2(q^d))dq/q\,, \,\omega_3=(\sum_{d|3} d \,f_3(q^d))dq/q\,, $$
 $$\omega_4=\frac{f_4+f_5}{2} dq/q\text{ and } \omega_5=\frac{f_5-f_4}{2 \sqrt 2} dq/q\,,
 $$
 we get  that $\cL_2=\{Q=a\, Q_1+b\,Q_2+c\,Q_3\colon a,b,c\in\Q\}$, where
 $$
 \begin{array}{cr}
 Q_1=&6 x_1^2 + 5 x_1 x_2 + 7 x_1 x_3 - 11 x_2 x_3 - 9 x_3^2 + 2 x_4^2 + 48 x_4 x_5 +
 16 x_5^2\,,\\
 Q_2=& 2 x_1 x_2 + 3 x_2^2 - 2 x_1 x_3 + 4 x_2 x_3 - 3 x_3^2 - 4 x_4^2 + 16 x_5^2 \,,\\
 Q_3=&x_2 x_4 - x_3 x_4 - 2 x_1 x_5 + x_2 x_5 + x_3 x_5\,.
 \end{array}
 $$
 Now, $\cL_2^{ns}=\langle Q_3\rangle$. For all possible choices $(i,j)\in\{(1,2),(1,3),(2,3),(4,5
 )\} $,   we only have the possibility $(i,j)=(4,5)$ satisfying condition (\ref{eq2}), i.e.
 $Q-Q(x_1,x_2,x_3,-x_4,-x_5)\in \cL_2^{ns}$ for all $Q\in\cL_2$.  Therefore, there is a unique non trivial involution $u$ sending $(\omega_1,\omega_2,\omega_3,\omega_4,\omega_5)$ to $\pm(\omega_1,\omega_2,\omega_3,-\omega_4,-\omega_5)$.  The number of eigenvalues equal to $1$ is the genus of $X_u=X_0^*(N)/\langle u\rangle$ and the corresponding eigenvectors among the differentials $\{\omega_i, i\leq 5\}$ is a basis of the pullback of $\Omega^1_{X_u}$.  To decide the precise  sign,
  we compute the number of fixed points of $u$. The set of such points in $\PP^4(\overline{\Q})$   is the set of points  of the form $(0, 0, 0, x_4, x_5) $ and of the form $ (x_1, x_2, x_3, 0, 0)$ satisfying $Q_i(x_1,\cdots,x_5)=0$ for $1\leq i\leq 3$. Since this set has  $4$ points, all of them of the form $x_4=x_5=0$,   the genus of $X_u$  is $2$ and $\omega_4$ and $\omega_5$ is a basis of the pullback of $\Omega^1_{X_u}$. Put
 $X=\omega_4/\omega_5$ and $Y=dx/\omega_5$ and we obtain
 $$Y^2=X^6 + 8 X^4 + 20 X^2 + 12 X + 4\,.$$

   For $g_N^*= 6$, we get $\dim \cL_2^{ns}=0$. The genus of a quotient curve by an involution  must be $2$  or $3$. For the possible choices of pairs $(i,j)$ or triples $(i,j,k)$, there are  polynomials in $\cL_2$ that are not  simultaneously even in the corresponding variables. Hence, all curves $X_0^*(N)$ of genus $6$ have trivial automorphism group.

 For $g_N^*\geq 6$, in all cases $\dim \cL_2^{ns}=(g-1)(g-6)/2$. We have to consider all choices $(i_1,\cdots,i_r)$ such that $ (g_N^*-5)/2\leq r\leq (g_N^*+1)/2 $ and $\omega_{i_1},\cdots,\omega_{i_r}$ is a basis of the puulack of  regular differentials of  a quotient of $J_0^*(N)$ of dimension $r$. After computing equations, we can claim that in all these cases there are polynomials in $\cL_2$ not satisfying the condition (\ref{eq2}).

 Nevertheless, we note that some of these computations can be simplified   by using Proposition  and Remark \ref{nopetri}. Indeed, for the pairs $(N,p)$ in the set
 $$\{(1155,3),(705,3),(987,3),(1365,3), (1365,5),(957,3),(1055,5)\}\,,$$
 the curve $X_0^*(N/p)$ has trivial automorphism group and  $X_0^*(N/p)/\F_p$ is not hyperelliptic. Hence, the choice $(i_1,\cdots,i_r)$ should
 contain the  variables corresponding to the basis of the pullback of  $\Omega^1_{J_0^*(N/p)}$. Next, in the four following examples, we show how we use these results.

 For instance,  if $X_0^*(1365)$ of genus $9$ has a nontrivial involution $u$, by Porposition\ref{restrict}, the jacobian of the quotient curve $X_u$ has  $J_0^*(273) $ and $J_0^*(455)$ as factors. This fact leads to  the contradiction that $\Jac(X_u)$ has a factor of dimension $6$ when the genus of $X_u$ is at most $5$.

 For $N=957$, take a basis $\{g_i\}_{1\leq i \leq 4} $  of $S_2^*(319)\cap\Q[[q]]$ and consider the set of the cusp forms $h_i(q)=g_i(q)+3g_i (q^3)\in S_2^*(957)$. By observing the splitting of $J_0^*(957)$, the vector space spanned by $h_i(q) dq/q$ would be the pullback of the regular differential of the unique possible quotient curve and the vector space of the homogenous polynomial in $\Q[x_1,\cdots,x_{4}]$ of degree $2$ vanishing at these differentials  would have dimension $1$. After a computation, the dimension obtained is $0$. So, we can exclude the level $957$. Note that for $X_0^*(957)$, $\dim \cL_2=36$.

 For $N=705$, take a basis $\{g_i\}_{1\leq i \leq 5} $  of $S_2^*(235)\cap\Q[[q]]$ and consider the set of the cusp forms $h_i(q)=g_i(q)+3g_i (q^3)\in S_2^*(957)$. The vector space spanned by $h_i(q) dq/q$ should be the pullback of the regular differential of the unique possible quotient curve and the vector space of the homogenous polynomial in $\Q[x_1,\cdots,x_{5}]$ of degree $2$ vanishing at these differentials  should have dimension $3$. After a computation, the dimension obtained is $0$. So, we can exclude the level $705$.

  If $X_0^*(1551)$ has a non trivial involution, then the jacobian of the genus six curve $X_0^*(211)$ is  isogenous to the jacobian of the quotient curve. Since $|X_0^*(1555)(\F_2 )|-2|X_0^*(211)(\F_2)|= 4>0$, we can discard $N=1551$.
\hfill $\Box$

\section{Even levels}
The next proposition is the key result to delimit  the  even levels.
\begin{prop}\label{negative} Let $N>1$  be an odd square-free integer. Then,
$$ g_{2 N}^*-2 g_N^*\leq 2 \,.$$
Moreover, $g_{2 N}^*-2g_N^* <-1 \text{ for all } N>1239$ and, for
the particular case  $g_N^*>2$, one has that
\begin{itemize}
\item[\rm{(i)}]  $g_{2 N}^*=2\,g_N^*-1$ if, and only if, $N$ is in the set
$$
\begin{array}{c}
\{109, 113, 139, 151, 203, 227, 259, 263, 319, 355, 411, 445, 451,
455, 461,\\ 491, 505, 521, 555, 573, 581, 591, 695, 699, 779, 1001,
1131, 1239\}\,.
\end{array}
$$
\item[\rm{(ii)}]  $g_{2 N}^*=2\,g_N^*$ if, and only if, $N$ is in the set
$$\begin{array}{c}
\{173, 267, 281, 295, 339, 341, 359, 371, 377, 413, 419, 429, 431, \\447,479, 483, 501, 551,623, 627, 645, 663, 671, 755, 789, 987\}\,.
\end{array}
$$
\item[\rm{(iii)}]  $g_{2 N}^*=2\,g_N^*+1$ if, and only if, $N$ is in the set
$$
\{149, 179, 239, 249, 251, 269, 305, 311, 321, 329, 393, 395,519,
545, 689, 861, 897\}\,.
$$
\item[\rm{(iv)}] $g_{2 N}^*=2\,g_N^*+2$ if, and only if, $N=303$.
\end{itemize}

\end{prop}
\noindent {\bf Proof.} Let us denote by $\cP$ the set of integer primes dividing  $N$ and set $n:=|\cP|$. We have that the genera of $X_0(N)$ and $X_0(2N)$ are
\begin{equation}\label{fg}
g_N=1+\frac{\psi(N) }{12}-\frac{\nu_2}{4}-\frac{\nu_3}{3}-2^{n-1}\,,\quad g_{2N}=1+\frac{\psi(N) }{4}-\frac{\nu_2}{4}-2^{n}\,,
\end{equation}
where
$$\nu_2=\left\{\begin{array}{cr} 0& \text{if $\exists p\in\cP\,, p\equiv-1 \pmod 4$,}\\2 ^n & \text{otherwise.} \end{array} \right. \,,\quad
\nu_3=\left\{\begin{array}{cr}
0& \text{if $\exists p\in\cP\,, p\equiv-1 \pmod 3$,}\\ 2^{n-v_3(N)} &\text{otherwise.} \end{array} \right. \,,
$$
where $v_3$ denotes the $3$-adic valuation.
The genera of $X_0^*(N)$ and $X_0^*(2N)$ are
\begin{equation}\label{genus}
g_N^*=1+\frac{1}{2^n}(g_N-1)-\frac{1}{2^{n+1}}\sum_{1< d|N}\nu(N,d)\,,\quad g_{2N}^*=1+\frac{1}{2^{n+1}}(g_{2N}-1)-\frac{1}{2^{n+2}}\sum_{1< d|2 N}\nu(2N,d)\,,
\end{equation}
where $\nu(M,d)$ denotes the number of fixed points of $X_0(M)$ by the Atkin-Lehner involution $w_d$. Hence  $g_{2N}^*-2g_N^*+1$ is equal to
$$
\frac{1}{2^{n+1}}\left(-\frac{\psi(N)}{12} +\frac{3\,\nu_2}{4}+\frac{4\,\nu_3}{3}
-\frac{1}{2}\nu(2N,2)+\sum_{1<d|N} 2\nu (N,d)-\frac{1}{2}\nu(2 N,d)-\frac{1}{2}\nu(2N,2d)+2^n\right)\,.
$$
Let $s=0,1$. By \cite{Kluit}, we know when $d\geq 5$:
$$
\nu(d,d)=\left\{ \begin{array}{lc} h(-4d)\,, & \text{if $d\not\equiv -1 \pmod 4$,}\\h(-4d)+ h(-d)\,, & \text{if $d\equiv -1 \pmod 4$,}
\end{array}  \right.
$$
$$
\nu(2 d,d)=\left\{ \begin{array}{lc} h(-4d)\,, & \text{if $d\not\equiv -1 \pmod 4$,}\\h(-4d)+ 3h(-d)\,, & \text{if $d\equiv -1 \pmod 4$,}
\end{array}  \right.
$$
where $h(D)$ is the class number of the order of  discriminant $D$ of a quadratic field, and for  $d|N$:
$$
\nu (2^s\, N,d)=\prod_{p|N/d}\left( 1+\left(\frac{-d}{p} \right) \right)\nu(2^s d, d)\,, \quad \nu (2N,2d)=\prod_{p|N/d}\left( 1+\left(\frac{-d}{p} \right) \right)\nu(2d,2d)\,.
$$
For $d<5$:
$$\nu(2N,2)=\prod_{p|N}\left( 1+\left(\frac{-1}{p} \right) \right)+\prod_{p|N}\left( 1+\left(\frac{-2}{p} \right) \right)\,,  \quad \nu(2^s3N,3)=2\prod_{p|N}\left( 1+\left(\frac{-3}{p} \right) \right) \,.
$$\,
We also know that for $d\equiv -1 \pmod 4$, $h(-4d)$ is $h(-d)$  or $3h(-d)$ depending on wether  $d\equiv -1\pmod 8$ or not (see \cite{Cox}[Theorem 7.24].

Since  $\nu_2\leq \nu(2N,2)$ and $\nu_2,\nu_3\leq 2^n$,  we have
$$
g_{2N}^*-2g_N^*+1\leq \frac{1}{2^{n+1}}\left(-\frac{\psi(N)}{12} +\sum_{1<d|N} (2\nu(N,d)-\frac{1}{2}\nu(2N,d))+\frac{31}{12} 2^n\right)\,.
$$

If $D$  is the discriminant of an order of an imaginary quadratic field, we known $h(D)\leq \frac{1}{\pi}|D|^{1/2}\log(|D|)$  (see Appendix in \cite{Serre}). One the one hand,
$$2\nu(d,d)- \frac{1}{2}\nu(2,d)=\left\{ \begin{array}{rr} \frac{3}{2} h(-4d)&\text{ if $d\not \equiv  1\pmod 4$,}\\
 2 h(-d)&\text{ if $ d \equiv 7\pmod 8$,}\\
 \frac{5}{3}h(-4d)&\text{ if $ d \equiv  3\pmod 8$,}\\
\end{array} \right.\,,$$
   and, on the other hand, $\log (4 d)\leq d^{1/4}+3$.
   Hence, we get
$$ 2\nu(d,d)- \frac{1}{2}\nu(2,d)\leq \frac{10}{3\pi} (d^{1/2}+3 d^{1/4})\,.$$  Since $\displaystyle{\frac{31}{12} < \frac{10}{3 \pi} (1^{3/4}+3\cdot  1^{1/2})}$, we have
$$
g_{2N}^*-2g_N^*+1<\frac{1}{2^{n+1}}\left(-\frac{\psi(N)}{12}+\frac{10}{3\pi}\sum_{d|N} 2^{\omega(N/d)}(d^{3/4}+ 3\,d^{1/2}) \right)\,.
$$
Since
$$
\sum_{d|N} 2^{\omega(N/d)}d^{3/4}=\prod_{p\in \cP}2+p^{3/4}\,,\quad \sum_{d|N} 2^{\omega(N/d)}d^{1/2}=\prod_{p\in \cP}2+p^{1/2 }\,,
$$
the following inequality
\begin{equation}\label{ineq}
\frac{10}{3\pi}\left (\prod_{p\in\cP}\frac{2+p^{3/4}}{1+p}+3\prod_{p\in\cP}\frac{2+p^{1/2}}{1+p}\right)-\frac{1}{12}<0
 \end{equation}
implies $g_{2N}^*-2g_N^*+1 <0$.
Write $\cP=\{p_1<\cdots<p_n\}$.  For two odd square-free integers $N=\prod_{i=1}^np_i$ and $N'=\prod_{i=1}^{n'}p_i'$, we define the order $N\preceq N'$ if, and only if,  $n\leq n'$ and $p_i\leq p_i'$ for al $i\leq n$. The real function $f(x)=\frac{2+x^\alpha}{1+x}$ with $ \frac{1}{2}\leq \alpha \leq \frac{3}{4}$ is decreasing for $x\geq3$ and $f(x)<1$ for $x\geq 5$. Hence, if the  the inequality (\ref{ineq}) is satisfied   for $N$, then it holds  for all integers $N'\succeq N$. Therefore, the inequality (\ref{ineq}) is right  for the values $N$  that satisfy some of the following conditions
 \begin{enumerate}
  \item  $\omega(N)\geq  8$.
  \item For  $\omega(N)= 7$ when   $p_1\geq 5$ or $p_7 > 73$.
   \item For  $\omega(N)=6$ when $p_1\geq 7$ or $p_6 > 569$.
  \item For  $\omega(N)=5$ when  $p_1\geq 13$ or $p_5> 3373$.
  \item For  $\omega(N)=4$ when $p_1\geq 23$ or $p_4 > 16573$.
  \item  For  $\omega(N)=3$ when  $p_1\geq 53$ or $p_3> 37993$.
  \item For  $\omega(N)=2$, when $p_1\geq 269$ or $p_2 > 63737$.
  \item  For  $\omega(N)=1$ when  $p_1 >54277$ .
\end{enumerate}
The statement is obtained after computing $g_{2N}^*-2g_N^*$ for the remaining   values of $N$. \hfill $\Box$

\begin{rem} The proof presented in the above proposition needs a laborious computation because it involves a lot of possibilities and, thus,  many class numbers. In opinion of the authors, it has to exist a deeper explanation to justify the  inequality stated for $g_{2N}^*-2\,g_N^*$.
\end{rem}

\subsection{Candidates $X_0^*(2N)$ with non trivial automorphism group}
Combining Propositions  \ref{restrict} and \ref{negative},  the cases $g_{2N}^*-2\,g_N^*<-1$ can be discarded, because of the inequality $g_N^*>(g_{2N}^*+1)/2$ (see Remark \ref{restgenus}).
Now, we only have to study the automorphism group  of $X_0^*(2N)$  for the odd values $N$ contained in the following five lists:
\begin{enumerate}

\item The list of values of  $N$ with  $g_N^*\leq 2$ such that $g_{2N}^*>3$:
$$\begin{array}{cc}\ell_1=&\{ 101, 107, 131, 161,  167, 177, 191,205, 209, 213, 221,  285, 287,299 ,357\}\end{array}$$
By Lemma \ref{criterion1} we discard $N=191$ at $p=5$. For the
remaining $N$, the splitting of the jacobian of $X_0^*(2N)$ is {{
$$
\begin{array}{c|c|}
N&  J_0^*(2N) \\\hline
101 &1_{101}+3_{202}\\[3 pt]
131&1_{131}+1_{262}+2_{262}\\[ 3pt]
107  &2_{107}+1_{214}+1_{214}\\[ 3pt]
161 & 2_{161}+2_{322}\\ \hline
\end{array}
\begin{array}{c|c||}
N&  J_0^*(2N) \\\hline
167 &2_{167}+2_{314}\\[3 pt]
177&1_{118}+2_{177}+1_{354}\\[ 3pt]
213  &1_{142}+2_{213}+1_{426}\\[ 3pt]
285 & 1_{57}+1_{190}+1_{285}+1_{570}\\ \hline
\end{array}$$
$$
\begin{array}{c|c|}
N&  J_0^*(2N) \\\hline
205&1_{82}+2_{205}+2_{410}\\[ 3pt]
209  &2_{209}+3_{418}\\[ 3pt]
221 & 2_{221}+3_{442}\\ \hline
\end{array}\begin{array}{c|c|}
N&  J_0^*(2N) \\\hline
287 &1_{82}+2_{287}+1_{574}+1_{574}\\[3 pt]
357&1_{102} +1_{238}+2_{357}+1_{714}\\[3 pt]
\hline\hline
299  &2_{299}+4_{598}\\\hline
\end{array}
$$}}
The genus of a quotient curve of any of these curves $X_0^*(2N)$ by an involution non bielliptic must be equal to $2$ when $g_{2N}^*=4$ and $2$ or $3$  when $5\leq g_{2N}^*\leq 6$. Hence, we can discard $N=101$ and for the remaining cases we apply Lemma \ref{con}. Finally, we get that all curves $X_0^*(2 N)$ with $N\in\ell_1$ have trivial automorphism group.
\item The list containing the values of $N$  with $g_N^*>2$ such that $X_0^*(N)/\F_2$  could be hyperelliptic:
$$\ell_2=\{183,185,187,203,335,345,385\}\,.$$
The splitting of the jacobian of $X_0^*(2N)$ {{
$$
\begin{array}{c|c|}
N&  J_0^*(2N) \\\hline
183 &1_{61}+ 1_{126}+2_{183}\\[4 pt]
185&1_{37}+1_{185}+1_{185}+1_{370}\\[ 3 pt]
187  &3_{187}+1_{374}\\[3 pt]
 &   \\\hline
\end{array}
\begin{array}{c|c|}
N&  J_0^*(2N) \\\hline
203 &1_{58}+ 3_{203}+1_{406}\\[ 3 pt]
385 &1_{77}+1_{154}+3_{385}\\ \hline\hline
335&2_{67}+2_{335}+2_{670}\\[ 3 pt]
345  &2_{115}+1_{138}+2_{345}+1_{690}\\ \hline
\end{array}$$}}
 The curves $X_0^*(2N)$ with $N\in\ell_2\backslash\{183\}$ have trivial automorphism group. Indeed, we can discard $N=187$ because $J_0^*(374)$ does not have any two dimensional quotients and the remaining values of $N$ can be excluded by applying Lemma \ref{con}. For $X_0^*(366)$ the automorphism group has order  $2$ (cf. \cite{BaGon}).

\item  The list containing the values of $N$  with $g_N^*>2$ such that $\Aut(X_0^*(N))$ in not trivial, except to $183$  because it is  in the list $\ell_2$:
$$\ell_3=\{ 249, 303, 455\}\,,$$
with
$$
\begin{array}{c|c|c|c|}
N& 249  &  303 &455\\\hline
J_0^*(2N) &1_{83}+1_{166}+1_{249}+1_{249}+3_{498}& 1_{101}+3_{202}+2_{303}+2_{606}&1_{65}+1_ {91}+1_{455}+2_{910}\\ \hline
\end{array}$$
By applying Lemma \ref{con}, we get that all these three curves $X_0^*(2N)$ have trivial automorphism group. We point out that $X_0^*(202)/\F_3$ is not hyperelliptic (see Remark \ref{hip}) and the automorphism group of $X_0^*(202)$ is trivial. Hence,  by Proposition  \ref{restrict}, for $X_0^*(606)$ we only have to consider as
 a possible quotient those that is jacobian correspond to the factor
$J_0^*(202)$.

\item   The list $\ell_4$ containing the values of $N$ with  $g_N^*>2$ and $-1\leq g_{2N}^*-2g_N^*\leq  0$, except to $N=203$ that is in the list $\ell_2$. So, the list $\ell_4$ is the union of the sets
$$
\begin{array}{c}
\{109, 113, 139, 151,  227, 259, 263, 319, 355, 411, 445, 451, 455,
461, 491,\\ 505, 521, 555, 573, 581, 591, 695, 699, 779, 1001, 1131,
1239\}\,,
\end{array}
$$
for which $g_{2N}^*-2g_N^*=-1$, and
$$\begin{array}{c}\{173, 267, 281, 295, 339, 341, 359, 371, 377, 413, 419, 429, 431, \\ 447,
479, 483, 501, 551,
623, 627, 645, 663, 671, 755, 789, 987
\}\,.\end{array}$$
For all these values, $X_0^*(N)$ has trivial  automorphism group and its reduction  modulo $2$ is not hyperelliptic (see Lemma \ref{hyp2}). By Proposition \ref{restrict}, if $X_0^*(2N)$ has a non trivial involution, then the jacobian of the quotient curve is isogenous to  $J_0^*(N)$ because $g_N^*$ is the greatest genus of a quotient curve of $X_0^*(2N)$.

{{ When $g_{2N}^*-2g_N^*=-1$, if $X_0^*(2 N)$ has a non trivial
involution, by Hurwitz, this cannot have fixed points. Hence, for
all odd primes $p\nmid N$ and all integers $k>0$, the number
$R(2N,p^k):=|X_0^*(2N)(\F_{p^k})|$ must be even. Since
$$
\begin{array}{c|c|c||}
N& p^k&R(2N,p^k)\\ \hline 109& 3^3& 21\\ 113&3 & 9\\ 139& 3^3& 21\\
151&3 &17 \\ 227&3 & 9\\ 259& 5&17\\  263&3^3 &17 \\ \hline
\end{array}
\begin{array}{c|c|c||}
N& p^k&R(2N,p^k)\\ \hline
 319& 3&9\\ 355& 3&7 \\ 411&5^5 &3323 \\ 445&3^3 &29 \\ 451& 3^3& 11 \\455&3 &7 \\ 461&3 &15\\\hline
\end{array}
\begin{array}{c|c|c||}
N& p^k&R(2N,p^k)\\ \hline
  491& 3^5&75 \\505& 3^5&147 \\ 521&3^3&29 \\ 555&11 & 17\\ 573 & 5^3&89\\ 581 &3^3 &27\\ 591 & 5&17\\\hline
\end{array}
\begin{array}{c|c|c||}
N& p^k&R(2N,p^k)\\ \hline
  695&3 &9 \\699& 5&23 \\ 779&3^3 &21 \\ 1001& 3^5&305 \\1131 & 5^7&75993\\ 1239&5 &15\\  & &\\\hline
\end{array}
$$
we can discard all these values.

When $g_{2N}^*=2g_N^*$, by applying Lemma \ref{criterion1}, we can discard  some values of $N$. With the prime $p=3$, the values of $N$ in the set
$$\begin{array}{c}\{  173,281, 359,377,419,431, 479,671,755\}\,,\end{array}$$
and with the prime $p=5$, the following values
$$\{413,501, 623,789\}\,.$$

}}






{{ We only have to consider
$$\begin{array}{c|c||}
N&  J_0^*(2N) \\\hline
295 &1_{118}+3_{295}+2_{590}\\[3 pt]
429&1_{143}+1 _{286}+2_{429}+2_{858} \\\hline\hline
267 &1_{89}+ 2_{178}+3_{267}+2_{534}\\[3 pt]
341&4_{341}+4_{682}\\[ 3 pt]
483  &2_{161}+2_{322}+2_{483}+1_{966}+1_{966}\\\hline\hline
339 &3_{113}+ 2_{226}+2_{339}+1_{678}+2_{678}\\[3 pt]
371&1_{53}+1_{106}+1_{371}+3_{371}\\ & +1_{742}+3_{742}\\\hline
\end{array}
\begin{array}{c|c|}
N&  J_0^*(2N) \\\hline
551 &1_{58}+ 2_{551}+ 3_{551}+2_{1102}+2_{1102}\\[3 pt]
645&1_{43}+1_{129}+1_{215}+1_{258}\\[3 pt]
 &+1_{430}+2_{645}+3_{1290}\\[ 3 pt]
663  &1_{102}+2_{221}+3_{442}+3_{663}+1_{1326}\\\hline\hline
447 &3_{149}+ 1_{298}+ 3_{298}+3_{447}+1_{894}+1_{894}\\[3 pt]
627&1_{57}+2_{209}+3_{418}+3_{627}+3_{1254}\\ \hline\hline
987  &1_{141}+2_{282}+3_{329}+4_{658}+2_{987} \\ & +3_{987}+3_{1974}\\\hline
\end{array}
$$}}

 {{

It can be checked that he curves  $X_0^*(322)/\F_3$, $X_0^*(226)/\F_3$,
$X_0^*(430)/\F_3$, $X_0^*(442)/\F_3$, $X_0^*(298)/\F_3$, $X_0^*(418)/\F_3$, and $X_0^*(658)/\F_3$ are not hyperelliptic. Since the pullback of the regular differentials of a quotient curve corresponds to the pullback of   $\Omega^1_{J_0^*(N)}$, by Proposition  \ref{restrict},
we can discard the values $$N\in\{483,339, 645 ,663,447,627,987\}\,. $$
Although $X_0^*(178)$ has a unique non trivial involution,   whose quotient curve  is isogenous to an elliptic curve of conductor $89$, the curve $X_0^*(534)$ can also be discarded. Indeed,  a non trivial involution of $X_0^*(534)$ should induce the identity on the elliptic curve of conductor $89$ and, thus on the curve $X_0^*(178)$. But this fact leads to the contradiction that the jacobian of the quotient  curve would have as a factor an abelian surface attached to a newform of level $178$.
By the remaining values, i.e.
$$N\in\{ 295,429,  341,371, 551\}\,,
$$
we apply Lemma \ref{con} and  we get that the five  curves $X_0^*(2N)$ have trivial automorphism group.}}


\item The list of values  $N$ with $g_N^*>2$ and $g_{2N}^*>2g_N^*$, except to  $N=249$ and $N=303$ because  they are   in the lists $\ell_3$, are the following:
$$
\ell_5=\{149, 179, 239,  251, 269, 303, 305, 311, 321, 329, 393, 395,
519, 545, 689, 861, 897\}\,.
$$
For all these cases, $g_{2N}^*-2g_N^*=1$. By Proposition
\ref{restrict}, if $X_0^*(2N)$ has a non trivial involution, then
the pullback of the regular differentials of the quotient curve is
the pullback of $\Omega^1_{J_0^*(N)}$ or $\Omega^1_{J_0^*(N)\times
E}$ for some elliptic curve  $E$ in the decomposition of
$J_0^*(2\,N)$, {{whose conductor does not divide $N$}}.

By applying Lemma \ref{criterion1} to $X_0^*(2\, N)$, we can discard the values of $N$ in the set
$$\{  149,179,239, 251, 269, 311, 393,519,545\}\,,$$
taking $p=5$ for $2 \,N=2\cdot 393, 2\cdot 519, 2\cdot 179$, $p=7$
for $2 \,N=2\cdot 545$ and $p=3$ for the remaining cases. The
splitting of $J_0^*(2N)$ for the non discarded cases  is as follows
{{
$$
\begin{array}{c|c||}
N&  J_0^*(2N) \\\hline
  329 &  3_{329}+ 4_{658}\\ \hline\hline
305 &1_{61}+ 1_{122}+3_{305}+4_{610}\\[3 pt]
321&2_{107}+1_{214}+1_{214}+2_{321}+3_{642}\\[ 3 pt]
395  &1_{79}+1_{158}+3_{395}+4_{790}\\ &\\\hline
 \end{array}
\begin{array}{c|c||}
N&  J_0^*(2N) \\\hline
861&1_{82}+1_{123}+1_{246}+2_{287}+1_{574}+1_{574}\\[3 pt] & +4_{861}+1_{1722}+1_{1722}+2_{1722}\\[3 pt]
897&1_{138}+2_{299}+4_{598}+5_{897}+3_{1794}\\\hline\hline
  689 & 1_{53}+1_{106}+1_{689}+2_{689}+3_{689}\\ &+3_{689}+ 3_{1378}+ 6_{1378}\\\hline
\end{array}
$$}}
It can be checked that the curves  {{$X_0^*(214)/\F_3$}},
$X_0^*(574)/\F_3$, and $X_0^*(598)/\F_3$ are not hyperelliptic. By
Proposition \ref{restrict}, the values $N\in\{321, 861, 897\}$
 can be discarded.

By the remaining values, i.e. $N=329,305,395, 689$, we apply Lemma
\ref{con} and  we get that the five   curves $X_0^*(2N)$ have
trivial automorphism group.

{{For instance, for $N=689$ we can proceed in a easier way. If
$X_0^*(1378)$ has a non trivial involution $u$, then  the jacobian
of the quotient curve $X_u$  is isogenous to $J_0^*(689)$ or
$J_0^*(689)\times A_g$,  where $g$ in the only newform in
$\New_{106}^*$.  None of these curves could be  hyperelliptic
because $|X_u(\F_3)|>8$. Take a basis $g_1,\cdots,g_9$ of
$S_2^*(689)\cap\Q[[q]]$. Put $h_i(q):=g_i(q)+2g_i(q^2)$ and
$h_{10}(q):=g(q)+13 g(q^{13})$. It can be checked that the dimension
of the vector space of homogenous polynomials in  $9$ variables
(resp. $10$ variables) vanishing at $h_1,\cdots,h_9$ (resp.
$h_1,\cdots,h_{10} $) has dimension $1$. Hence, $N=689$ can be
excluded. Note that for $X_0^*(1378)$, $\dim \cL_2=153$.}}
\end{enumerate}

As a consequence of this analysis  and taking into account  \cite{BaGon}[Proposition 24 ], we get  the following result.
\begin{prop}\label{even}
Let $N>1$ be an {{even}} square-free  integer such that $g_N^*>3$
and $X_0^*(N)$ is not bielliptic. The group $\Aut (X_0^*(N))$ is non
trivial if, and only if, $N=366$. In this case, the order of this
group is $2$ and an equation  for the quotient curve by the non
trivial involution is given by
$$
Y^2 =X^6 - 6 X^5 + 23 X^4 - 42 X^3 + 53 X^2 - 24 X + 4\,.
$$
\end{prop}
This proposition together Proposition \ref{odd} conclude the proof of Theorem \ref{main2}.

\section*{Acknowledgements}
 We thank  C. Ritzenthaler for his   comments and  B.   Poonen,  whose  contribution was  decisive  to prove Proposition \ref{negative}.


\noindent{Francesc Bars Cortina}\\
{Departament Matem\`atiques, Edif. C, Universitat Aut\`onoma de Barcelona\\
08193 Bellaterra, Catalonia}\\
{francesc@mat.uab.cat}

\vspace{1cm}

\noindent
{Josep Gonz\'alez Rovira}\\
{Departament de Matem\`atiques, Universitat Polit\`ecnica de Catalunya EPSEVG,\\
Avinguda V\'{\i}ctor Balaguer 1, 08800 Vilanova i la Geltr\'u,
Catalonia}\\
{josep.gonzalez@upc.edu}

\end{document}